\documentclass[12pt,letterpaper,oneside]{article}

\usepackage{amsmath,amssymb,amsthm,mathrsfs}
\usepackage{cite}

\newtheorem{lemma}{Lemma}[section]

\newtheorem{theorem}[lemma]{Theorem}
\newtheorem{proposition}[lemma]{Proposition}

\newtheorem{definition}[lemma]{Definition}
\newtheorem{remark}[lemma]{Remark}

\newcommand{\dd}{\,\mathrm{d}}

\newenvironment{ack}
  {\par\noindent\textbf{Acknowledgments}.\enspace}
  {\par}

\title{Diameter Growth of $\mathrm{U}(2^N)$ under Cliff Metrics}
\author{%
  Xiaobo Liu,\,\,\, Lei Zheng%
}

\date{}

\begin{document}

\maketitle

\begin{abstract}
	Cliff metrics are right invariant Riemannian metrics on the unitary group $\mathrm{U}(2^N)$ introduced by Nielsen and his collaborators in the geometric approach for gate complexity in quantum computations.   
In this paper, we compute the exponential growth rate for the diameter of $\mathrm{U}(2^N)$ equipped with cliff metrics as $N \rightarrow \infty$. We prove that if the norm of most generalized Pauli matrices are sufficiently large, this growth rate coincides with the growth rate of maximal gate complexity obtained by Knill. 
\end{abstract}

\section{Introduction}
	
	A geometric approach to complexity in quantum computations was introduced by Nielsen in \cite{Nielsen2006Geometric}. This approach has been extensively studied in many articles (see, for example, \cite{AuzziEtAl2021, Brandt2009, BrownSusskind2018, BrownSusskind2019, BrownSusskindZhao2017, BulchandaniSondhi2021, CrapsEtAl2022, DowlingNielsen2008, GuDohertyNielsen2008, Lin2019, NielsenEtAl2006Geometry, Nielsen2006OptimalControl, RibeiroTrancanelli2026, Wu2022}). The basic idea for the geometric approach is to equip the unitary group $\mathrm{U}(2^N)$ with a right-invariant Riemannian metric and define the geometric complexity of a unitary operator as its distance from the identity matrix $\mathrm{Id}$. The diameter of $\mathrm{U}(2^N)$ under such a metric corresponds to the maximal geometric complexity of unitary operators. A distinguished class of such metrics are the so called cliff metrics. It was proven in  \cite{NielsenEtAl2006Geometry} that a polynomial of the geometric complexity defined using cliff metrics  is equivalent to the approximate gate complexity in quantum computations. In this paper, we compute the exponential growth rate for the diameter of $\mathrm{U}(2^N)$ under the cliff metrics as $N \rightarrow \infty$.

A right-invariant Riemannian metric on $\mathrm{U}(2^N)$ is determined by an inner product on its Lie algebra $\mathfrak{u}(2^N)$. A basis of $\mathfrak{u}(2^N)$ is given by the {\it generalized Pauli matrices}
\[
\tilde\sigma_{I} = -\sqrt{-1}\sigma_{i_1}\otimes\cdots\otimes\sigma_{i_N}, 
\] 
where $ I =(i_1, \ldots, i_N)$ with $i_j \in \{0,1,2,3\}$ for $j=1, \ldots, N$, $\sigma_0$ is the $2 \times 2$ identity matrix, and $\sigma_1, \sigma_2, \sigma_3$ are the standard Pauli matrices (defined by equation \eqref{eq:Pauli}). The \emph{weight} (or \emph{locality}) of $\tilde \sigma_I$ is the number of $i_j \neq 0$ in $I$. Given a real number $q > 1$, the {\it cliff metric} $g_q$ is the right invariant metric on $\mathrm{U}(2^N)$ such that $g_q(\tilde\sigma_{I}, \tilde\sigma_{J})=0$ if $I \neq J$, $g_q(\tilde\sigma_{I}, \tilde\sigma_{I})=q$ if $w(\tilde \sigma_{I}) >2$, and $g_q(\tilde\sigma_{I}, \tilde\sigma_{I})=1$ if $w(\tilde \sigma_{I}) \leq 2$. The number $q$ is called the {\it penalty} for the cliff metric $g_q$, which is usually chosen to be $q=e^{\lambda N}$ for some $\lambda > 0$ (see, for example, \cite{NielsenEtAl2006Geometry}). The distance function induced by the Riemannian metric $g_q$ is denoted by $d_q$. As $q \rightarrow \infty$, $g_q$ tends to a sub-Riemannian metric on $\mathrm{U}(2^N)$, called the {\it infinite cliff metric}, whose distance function is denoted by $d_\infty$.

There have been many interests in estimating the growth rate for the diameter of cliff metrics. In \cite{Brown2023LowerBound} and \cite{Brown2024PolynomialEquivalence}, Brown proved that the diameter  of $g_q$ is bounded below by $\min\{\sqrt{q},\, 2^N\}$ up to factors sub-exponential in $N$. He also conjectured in \cite{Brown2023LowerBound} that the diameter of the infinite-cliff metric grows as $2^{N + o(N)}$. In \cite{Fang2026VariableCliff}, Fang disproved this conjecture by establishing a lower bound of $2^{\frac{4}{3}N - o(N)}$. For $q = 2^{\lambda N}$ with $0 < \lambda < \frac{8}{3}$, Fang proved that the diameter of $g_q$ on the projective unitary group $\mathrm{PU}(2^N)$ is asymptotic to $\frac{\pi}{\sqrt{3}}\sqrt{q}$. In this paper, we show that the diameter of the infinite cliff metric grows as $2^{2N + o(N)}$. For  cliff metrics with $q = 2^{\lambda N}$, we show that the diameter grows as $2^{\min\{\lambda/2,2\}N + o(N)}$ as $ N \rightarrow \infty$. More precisely, we will prove the following:	
\begin{theorem}\label{thm:main-theorem-intro}
	For each integer $N \geq 3$, let $\mu_{2^N}$ be the normalized Haar measure on $\mathrm{U}(2^N)$.
	\begin{enumerate}
	\item[(1)] Fix $\lambda > 0$, and set $q_{N,\lambda} = 2^{\lambda N}$. Then
	\[
	\lim_{N \to \infty} \frac{1}{N}\log_2 \mathrm{diam}(\mathrm{U}(2^N), d_{q_{N,\lambda}}) = \min\{\lambda/2,2\}.
	\]
	Moreover, for every $\eta < \min\{\frac{\lambda}{2},2\}$, we have
	\[
\mu_{2^N}(\{x \in \mathrm{U}(2^N) \mid d_{q_{N,\lambda}}(\mathrm{Id}, x) < 2^{\eta N}\}) \longrightarrow 0, \qquad \text{as }N \to \infty.
	\]
	\item[(2)] For the infinite-cliff metric, we have
	\begin{equation} \label{eq:rate-infty}
	\lim_{N \to \infty} \frac{1}{N}\log_2 \mathrm{diam}(\mathrm{U}(2^N), d_\infty) = 2,
	\end{equation}
	and, for every $\eta < 2$,
	\[
	\mu_{2^N}(\{x \in \mathrm{U}(2^N) \mid d_{\infty}(\mathrm{Id}, x) < 2^{\eta N}\}) \longrightarrow 0, \qquad \text{as }N \to \infty.
	\]
	\end{enumerate}
\end{theorem}

		Theorem~\ref{thm:main-theorem-intro} allows a direct comparison between geometric complexity and quantum gate complexity. By classical results of Knill~\cite{Knill1995} (recalled in Section~2.3), both the maximal exact gate complexity and the maximal approximate gate complexity (for any fixed error $\epsilon \in (0,\sqrt{2})$) on $N$ qubits grow as $2^{2N + o(N)}$. Thus, for the infinite-cliff metric and for finite-cliff metrics with $q = 2^{\lambda N}$ ($\lambda \geq 4$), the diameter of $\mathrm{U}(2^N)$ has the same exponential growth rate as the maximal exact gate complexity.
		
		A further motivation comes from the conjectures of Brown, Freedman, Lin, and Susskind~\cite{BrownFreedmanLinSusskind2023}. Their Conjecture~1 predicts that within certain families of complexity metrics, there is a critical metric such that whenever the penalties are larger than those of the critical metric, the metrics have approximately the same distance function at large distances. Approximate agreement of their diameters is a necessary consistency check for this conjecture. Theorem~\ref{thm:main-theorem-intro} gives supporting evidence for this conjecture and
indicates that for the family of cliff metrics, the penalty of the critical metric is very likely given by
 $q= 2^{4 N}$.    

	The idea of obtaining a lower bound on diameter by estimating the Haar measure of cliff-metric balls was inspired by Fang's work~\cite{Fang2026SharpTypical, Fang2026VariableCliff}. The growth rate in equation \eqref{eq:rate-infty} was also conjectured for the infinite cliff metric on $\mathrm{PU}(2^N)$ by Fang in a  seminar talk.  Fang bounds the measure of a cliff-metric ball via  comparisons of Jacobian determinants. This is different from our method which uses the Dyson expansion to estimate the covering number of a cliff-metric ball by small balls in the standard bi-invariant metric. 

	 This paper is organized as follows. In Section~2, we collect the necessary preliminaries and establish diameter upper bounds. In Section~3, we derive covering estimates for finite-cliff balls. Finally, in Section~4, we combine these estimates to complete the proof of Theorem~\ref{thm:main-theorem-intro}.

\section{Preliminaries}

\subsection{Matrix norms and bi-invariant geometry on $\mathrm{U}(D)$}

	For an integer $D \geq 2$, let $M_D(\mathbb{C})$ denote the space of $D \times D$ complex matrices, and let $(\cdot)^\dagger$ denote the conjugate transpose. For a vertical vector $v \in \mathbb{C}^D$, its Euclidean norm is given by $|v|:=\sqrt{v^\dagger v}$. For $A \in M_D(\mathbb{C})$, its \emph{operator norm} is defined by
\[
\Vert A \Vert_{\mathrm{op}}:= \sup \bigl\{ |Av| \mid v \in \mathbb{C}^D, |v| = 1 \bigr\},
\]
and its \emph{Frobenius norm} (also called the \emph{Schatten $2$-norm}) is defined by
\[
\Vert A \Vert_{F} := \sqrt{\mathrm{tr}(A^\dagger A)}.
\]
We collect below the basic properties of the operator norm that will be used throughout the paper.

\begin{lemma}\label{lem:basic_property_of_operator_norm}
	Let $D \geq 2$ be an integer,  $U,V \in \mathrm{U}(D)$ and $A,B \in M_D(\mathbb{C})$. Then:
\begin{enumerate}
\item[(1)] $\Vert UAV \Vert_{\mathrm{op}} = \Vert A \Vert_{\mathrm{op}}$. In particular, $\Vert U \Vert_{\mathrm{op}} = 1$.

\item[(2)] $\Vert AB \Vert_{\mathrm{op}} \leq \Vert A \Vert_{\mathrm{op}} \Vert B \Vert_{\mathrm{op}}$ and $\Vert A + B \Vert_{\mathrm{op}} \leq \Vert A \Vert_{\mathrm{op}} + \Vert B \Vert_{\mathrm{op}}$.

\item[(3)] $D^{-1/2}\Vert A \Vert_F \leq \Vert A \Vert_{\mathrm{op}}$.

\item[(4)] For any $a \leq b$ and measurable map $C: [a,b] \to M_D(\mathbb{C})$ with $\int_a^b \Vert C(t) \Vert_{\mathrm{op}} \dd t < \infty$,
\[
\left\Vert \int_a^b C(t) \dd t \right\Vert_{\mathrm{op}} \leq \int_a^b \Vert C(t) \Vert_{\mathrm{op}} \dd t.
\]
\end{enumerate}
\end{lemma}

\begin{proof}
	Properties~(1)--(3) are standard properties of the operator and Frobenius norms; see \cite[Section~5.6]{HornJohnson2013}. To prove property~(4), for any $v \in \mathbb{C}^D$ with $|v| = 1$, 
we set $w := \bigl(\int_a^b C(t) \dd t \bigr) v$. If $w \neq 0$, setting $u:= \frac{w}{|w|}$, we obtain
\[
|w| = |u^\dagger w| = \left| \int_a^b u^\dagger C(t)v \dd t \right| \leq \int_a^b |u^\dagger C(t) v| \dd t \leq \int_a^b |u| |C(t)v| \dd t \leq \int_a^b \Vert C(t) \Vert_{\mathrm{op}} \dd t.
\]
The same bound holds trivially when $w = 0$. Taking the supremum over all $v \in \mathbb{C}^D$ with $|v| = 1$ yields~(4).
\end{proof}

	For each integer $D \geq 2$, let $\mu_D$ denote the normalized Haar measure on $\mathrm{U}(D)$. We equip $\mathrm{U}(D)$ with a bi-invariant Riemannian metric $g_1$ whose value at the identity is
\[
g_1(X,Y) := \frac{1}{D}\mathrm{tr}(X^\dagger Y),\qquad X, Y \in \mathfrak{u}(D).
\]
We denote the induced distance by $d_1$ and write $\Vert X \Vert_{g_1} = \sqrt{g_1(X,X)} = D^{-1/2}\Vert X \Vert_{F}$. 

	For $x \in \mathrm{U}(D)$ and $L > 0$, we write
\[
B_1(x, L) := \{y \in \mathrm{U}(D) \mid d_1(x,y) < L\}.
\]
For $A \subseteq \mathrm{U}(D)$ and $\epsilon > 0$, the $d_1$-covering number $N_1(A,\epsilon)$ is defined by
\[
N_1(A, \epsilon) := \min\Bigl\{m \in \mathbb{Z}_{\ge 0} \,\big|\, \exists\, y_1, \ldots, y_m \in \mathrm{U}(D),\, \bigcup_{i=1}^m B_1(y_i, \epsilon) \supseteq A\Bigr\}.
\]

	For $x,y \in \mathrm{U}(D)$, let $e^{\sqrt{-1}\theta_1}, \ldots, e^{\sqrt{-1}\theta_D}$ with $\theta_j \in (-\pi,\pi]$ be the eigenvalues of $x^\dagger y$. Rescaling the formula in \cite[Appendix~B.3, equation~(B.8)]{AubrunSzarek2017} by $D^{-1/2}$ gives
\begin{equation}\label{equ:geodesic_distance_bi_invariant}
d_1(x,y) = \left(\frac{1}{D}\sum_{j = 1}^D \theta_j^2 \right)^{1/2}.
\end{equation}
Since $|\theta_j| \leq \pi$ for every $j$, with equality attained when $y = -x$, it follows that
\begin{equation}\label{equ:diameter_of_bi_invariant_metric}
\mathrm{diam}(\mathrm{U}(D), d_1) = \pi.
\end{equation}

\begin{lemma}\label{lem:compare-distance}
	For every integer $D \geq 2$ and all $x, y \in \mathrm{U}(D)$,
\begin{equation}\label{equ:compare-distance}
D^{-1/2}\Vert x - y \Vert_F \leq d_1(x,y) \leq \frac{\pi}{2} D^{-1/2} \Vert x - y \Vert_F \leq \frac{\pi}{2} \Vert x - y \Vert_{\mathrm{op}}.
\end{equation}
\end{lemma}

\begin{proof}
	Rescaling the inequalities in \cite[Appendix~B.3, Exercise~B.5]{AubrunSzarek2017} by $D^{-1/2}$ gives
\[
D^{-1/2}\Vert x - y\Vert_F\leq d_1(x,y) \leq \frac{\pi}{2}D^{-1/2}\Vert x - y\Vert_F.
\]
Applying Lemma~\ref{lem:basic_property_of_operator_norm} yields the last inequality.
\end{proof}

\begin{lemma}\label{lem:unitary-bi-inv-riemannian-small-ball}
	There exist constants $\epsilon_0, c_0 > 0$ such that 
\begin{equation}\label{equ:unitary-bi-inv-riemannian-small-ball}
\mu_D(B_1(x,\epsilon_0)) \leq \exp(-c_{0}D^2)
\end{equation}
holds for every integer $D \geq 2$ and every $x \in \mathrm{U}(D)$.
\end{lemma}

\begin{proof}
Fix an integer $D \geq 2$. The space $(\mathrm{U}(D), d_{1})$ is compact and hence totally bounded. Thus, for each $\epsilon > 0$, we can choose a finite set $\Gamma = \{y_1, \ldots, y_s\}\subseteq \mathrm{U}(D)$ such that
\begin{equation}\label{equ:separated-condition}
\begin{aligned}
&d_{1}(y_i,y_j)>2\epsilon
&&\text{for all }i\neq j,\\
&\min_{1\leq i\leq s}d_{1}(x,y_i)\leq 2\epsilon
&&\text{for every }x\in\mathrm{U}(D).
\end{aligned}
\end{equation}
Consequently,
\[
\mathrm{U}(D) = \bigcup_{i = 1}^{s} B_{1}(y_i,3\epsilon).
\]

Claim: There is a universal constant $c > 0$ (independent of $D, \delta$) such that
\begin{equation}\label{equ:covering-number-lower-bound}
N_{1}(\mathrm{U}(D),\delta) \geq \left(\frac{c}{\delta}\right)^{D^2} \qquad \text{for every }\delta \in (0, \pi].
\end{equation}
Here $\pi = \mathrm{diam}(\mathrm{U}(D), d_1)$ by~\eqref{equ:diameter_of_bi_invariant_metric}.

For $x \in \mathrm{U}(D)$ and $r > 0$, set $\overline{B}_F(x,r):= \{y \in \mathrm{U}(D) \mid \Vert x - y\Vert_F \leq r\}$, and define the corresponding covering number by
\[
\overline{N}_F(\mathrm{U}(D), r):= \min\Bigl\{m \in \mathbb{Z}_{\geq 0}\,\big|\, \exists z_1, \ldots, z_m \in \mathrm{U}(D),\, \bigcup_{i = 1}^m \overline{B}_{F}(z_i,r) \supseteq \mathrm{U}(D)\Bigr\}.
\]
Applying \cite[Theorem~5.11 and Table~5.1]{AubrunSzarek2017} to $\mathrm{U}(D)$ equipped with the extrinsic distance induced by the Frobenius norm $\Vert \cdot \Vert_{F}$, we obtain a universal constant $c_1 > 0$ (independent of $D$ and $r$) such that
\[
\overline{N}_F(\mathrm{U}(D), r) \geq \left(\frac{2c_1\sqrt{D}}{r}\right)^{D^2}\qquad \text{for every }r \in (0,2\sqrt{D}].
\]

	By Lemma~\ref{lem:compare-distance}, we have
\[
D^{-1/2}\Vert x - y\Vert_F \leq d_1(x,y) \qquad \text{for all }x,y \in \mathrm{U}(D),
\]
and hence, for every $x \in \mathrm{U}(D)$ and $\delta > 0$,
\[
B_1(x, \delta) \subseteq \overline{B}_F(x,\sqrt{D}\delta).
\]
Consequently, for every $\delta \in (0,2]$, we have $\sqrt{D}\delta \in (0,2\sqrt{D}]$, and therefore
\[
N_1(\mathrm{U}(D), \delta) \geq \overline{N}_F(\mathrm{U}(D), \sqrt{D}\delta) \geq \left(\frac{2c_1\sqrt{D}}{\sqrt{D}\delta}\right)^{D^2} = \left(\frac{2c_1}{\delta}\right)^{D^2}.
\]
Setting $c := \min\{2c_1,2\}$, we obtain \eqref{equ:covering-number-lower-bound} for every $\delta \in (0,2]$. For $\delta \in (2, \pi]$, we have \linebreak$\frac{c}{\delta} < 1$, while $N_1(\mathrm{U}(D), \delta) \geq 1$ trivially, so \eqref{equ:covering-number-lower-bound} also holds in this case. This proves the claim.

	Applying \eqref{equ:covering-number-lower-bound} with $\delta = 3\epsilon$, whenever $0 < \epsilon \leq \frac{\pi}{3}$, we obtain
\[
|\Gamma| \geq N_{1}(\mathrm{U}(D),3\epsilon) \geq \left(\frac{c}{3\epsilon}\right)^{D^2}.
\]

Moreover, the balls $B_{1}(y_i, \epsilon)$ are pairwise disjoint.

By the left invariance of both the metric $g_1$ and the Haar measure,
\[
1 \geq \mu_D\left(\bigcup_{i = 1}^s B_{1}(y_i, \epsilon)\right) = |\Gamma|\,\mu_D(B_{1}(\mathrm{Id},\epsilon)).
\]
Combining these estimates gives
\[
\mu_D(B_{1}(\mathrm{Id},\epsilon)) \leq \frac{1}{|\Gamma|} \leq \left(\frac{3\epsilon}{c}\right)^{D^2}.
\]
By left invariance, the same estimate holds with $\mathrm{Id}$ replaced by any $x \in \mathrm{U}(D)$.

Finally, choose $\epsilon_0 > 0$ such that $3\epsilon_0 < c$ and $3\epsilon_0 \leq \pi$. Set
\[
c_0:= \ln \left(\frac{c}{3\epsilon_0}\right) > 0.
\]
Then
\[
\mu_D(B_{1}(x,\epsilon_0)) \leq \exp(-c_0 D^2),
\]
for every $x \in \mathrm{U}(D)$. This completes the proof of the lemma.
\end{proof}

\subsection{Cliff metrics and infinite-cliff geometry}

	Let $\sigma_0 = \begin{bmatrix}1 & 0\\0 & 1\end{bmatrix}$, and let
\begin{equation} \label{eq:Pauli}
\sigma_1 = \begin{bmatrix}0 & 1\\1 & 0 \end{bmatrix},
\qquad \sigma_2 = \begin{bmatrix}0 & -\sqrt{-1}\\ \sqrt{-1} & 0 \end{bmatrix},
\qquad \sigma_3 = \begin{bmatrix}1 & 0\\0 & -1 \end{bmatrix}
\end{equation}
be the standard Pauli matrices. Fix an integer $N \geq 3$ and set $D = 2^N$. Using the notation introduced in Section~1, we set
\[
\mathcal{I}_N := \{(i_1,\ldots, i_N) \mid i_1, \ldots, i_N \in \{0,1,2,3\}\},\qquad \mathcal{P}_N := \{I \in \mathcal{I}_N \mid w(\tilde \sigma_I) \leq 2\}.
\]
Since each $\sigma_j$ ($j \in \{0,1,2,3\}$) is unitary, each generalized Pauli matrix $\tilde\sigma_I = -\sqrt{-1}\sigma_{i_1} \otimes\cdots\otimes\sigma_{i_N}$ belongs to $\mathrm{U}(D)$, and hence by Lemma~\ref{lem:basic_property_of_operator_norm},
\begin{equation}\label{equ:unitary_Pauli_matrices}
\Vert \tilde\sigma_I \Vert_{\mathrm{op}} = 1,\qquad I \in \mathcal{I}_N.
\end{equation}
Moreover, $\{\tilde\sigma_I\}_{I \in \mathcal{I}_N}$ is a $g_1$-orthonormal basis of $\mathfrak{u}(D)$. Define the horizontal subspace
\[
\mathcal{H}_N := \mathrm{span}_{\mathbb{R}}\{\tilde \sigma_I \mid I \in \mathcal{P}_N\}.\]
The dimension of $\mathcal{H}_N$ is given by
\[b_N := 1 + 3N + 9\binom{N}{2},
\]
and $\{\tilde\sigma_I\}_{I \in \mathcal{P}_N}$ is a $g_1$-orthonormal basis of $\mathcal{H}_N$. Let $P_N: \mathfrak{u}(D) \to \mathcal{H}_N$ denote the $g_1$-orthogonal projection.

We now give the definition of the finite- and infinite-cliff metrics.

\begin{definition}\label{cliff-metric}
	With the notation above, for each $q \geq 1$, the \emph{cliff metric} $g_q$ introduced in Section~1 is the right-invariant Riemannian metric on $\mathrm{U}(D)$ given at the identity by
\[
g_q(X,Y) := g_1(X, P_N(Y) + q(Y - P_N(Y))), \qquad \textnormal{for all }X, Y \in \mathfrak{u}(D).
\]
(Note that setting $q = 1$ recovers the bi-invariant metric $g_1$.) We write $\Vert Xx \Vert_{g_q} = \sqrt{g_q(Xx,Xx)}$ for $X \in \mathfrak{u}(D), x \in \mathrm{U}(D)$, and denote the $g_q$-length of a $d_1$-Lipschitz curve $\gamma$ on $\mathrm{U}(D)$ by $l_{g_q}(\gamma)$.

Right translation of $\mathcal{H}_N$ gives a distribution $\Delta_{\mathcal{H}_N}$ on $\mathrm{U}(D)$ defined by
\[
\Delta_{\mathcal{H}_N}(x) = (R_x)_{*\mathrm{Id}}\mathcal{H}_N = \{Xx\mid X \in \mathcal{H}_N\} 
\subset T_x \mathrm{U}(D), \qquad x \in \mathrm{U}(D),
\] 
where $R_x(y) = yx$. $\Delta_{\mathcal{H}_N}$ is equipped with the right-invariant inner product $g_{\infty}(Xx, Yx) := g_1(X,Y)$ for $X,Y \in \mathcal{H}_N$. We write $\Vert v \Vert_{g_\infty} = \sqrt{g_{\infty}(v, v)}$ for $v \in \Delta_{\mathcal{H}_N}(x)$. Since $\mathcal{H}_N$ satisfies the Lie bracket-generating condition, the triple $(\mathrm{U}(D), \Delta_{\mathcal{H}_N}, g_\infty)$ is a sub-Riemannian manifold.
The sub-Riemannian metric $g_\infty$ is called the \emph{infinite-cliff metric}.
 A curve $\gamma: [0,1] \to \mathrm{U}(D)$, Lipschitz with respect to $d_1$, is called \emph{admissible} if there exists a control $u = (u_I)_{I \in \mathcal{P}_N} \in L^\infty([0,1]; \mathbb{R}^{b_N})$ such that
\[
\dot \gamma(t) = \left(\sum_{I \in \mathcal{P}_N} u_I(t)\tilde\sigma_I\right)\gamma(t), \qquad \textnormal{for almost every }t \in [0,1].
\]
Its length $l_{\infty}(\gamma)$ is defined by
\[
l_\infty(\gamma) := \int_0^1 \Vert \dot \gamma(t) \Vert_{g_\infty} \dd t = \int_0^1 \left(\sum_{I \in \mathcal{P}_N}u_I(t)^2\right)^{1/2}\dd t.
\]
The Carnot--Carath\'eodory distance is given by
\[
d_\infty(x,y) := \inf\{l_{\infty}(\gamma) \mid \gamma: [0,1] \to \mathrm{U}(D) \text{ is admissible},\, \gamma(0) = x,\, \gamma(1) = y\}.
\]
\end{definition}

\begin{remark}
We follow the convention of \cite[Chapter~3]{AgrachevBarilariBoscain2020}, using Lipschitz admissible curves with essentially bounded controls. Alternatively, one may allow absolutely continuous curves with controls in $L^1([0,1];\mathbb{R}^{b_N})$. Although this enlarges the class of admissible curves, it yields the same distance; see \cite[Proposition~3.72]{AgrachevBarilariBoscain2020}.
\end{remark}

	For $x \in \mathrm{U}(D)$ and $L > 0$, we also use the notation
\[
\begin{aligned}
B_q(x, L) &:= \{y \in \mathrm{U}(D) \mid d_q(x, y) < L\},\qquad q > 1,\\
B_\infty(x, L) &:= \{y \in \mathrm{U}(D) \mid d_\infty(x, y) < L\}.
\end{aligned}
\]
Since $g_1(X,X) \leq g_q(X,X) \leq q \, g_1(X,X)$ for all $X \in \mathfrak{u}(D)$, and 
for every admissible curve $\gamma$,  
$\Vert \dot\gamma(t) \Vert_{g_q} = \Vert \dot \gamma(t) \Vert_{g_\infty}$ $\textnormal{for almost every } t \in [0,1]$, we obtain
\begin{equation}\label{equ:metric_comparison}
d_1(x,y) \leq d_q(x,y) \leq \sqrt{q} \, d_1(x,y), \qquad d_q(x,y) \leq d_\infty(x,y),
\end{equation}
for all $x,y \in \mathrm{U}(D)$ and $q > 1$. Combining~\eqref{equ:diameter_of_bi_invariant_metric} with~\eqref{equ:metric_comparison} yields
\begin{equation}\label{equ:diameter_upper_bound_1_of_cliff_metric}
\mathrm{diam}(\mathrm{U}(D), d_q) \leq \pi \sqrt{q}.
\end{equation}

\subsection{Gate complexity and diameter upper bounds}

	Fix an integer $N \geq 3$ and set $D = 2^N$. We identify $\mathbb{C}^{D}$ with the $N$-qubit Hilbert space $(\mathbb{C}^2)^{\otimes N}$. To label the individual tensor factors, write $(\mathbb{C}^2)^{\otimes N} = \bigotimes_{j = 1}^N \mathcal{V}_j$, where $\mathcal{V}_j$ denotes the $j$-th copy of $\mathbb{C}^2$. For any subset $S \subseteq \{1, \ldots, N\}$ with complement $\overline{S} :=\{1, \ldots, N\} \setminus S$, define
\[
\mathcal{V}_S := \bigotimes_{j \in S} \mathcal{V}_j, \qquad \mathcal{V}_{\overline{S}} :=\bigotimes_{j \in \overline{S}} \mathcal{V}_j.
\]
Here the empty tensor product is understood as $\mathbb{C}$. Reordering the tensor factors gives a canonical unitary identification $(\mathbb{C}^2)^{\otimes N} \cong \mathcal{V}_S \otimes \mathcal{V}_{\overline{S}}$. Following~\cite[Section~2]{Knill1995} (see also~\cite[Chapter~4]{NielsenChuang2010}), for an integer $1 \leq k \leq N$, a unitary operator $U \in \mathrm{U}(D)$ is said to have \emph{locality at most $k$} if there exist a subset $S \subseteq \{1, \ldots, N\}$ with $|S| = k$ and an operator $W \in \mathrm{U}(\mathcal{V}_S)$ such that 
$
U = W \otimes \mathrm{Id}_{\mathcal{V}_{\overline{S}}}
$
after reordering.
Let $E \subseteq \mathrm{U}(D)$ denote the set of unitary operators of locality at most $2$.

	Following~\cite[Section~II.A]{Nielsen2006OptimalControl}, for $U \in \mathrm{U}(D)$, the \emph{exact gate complexity} of $U$ is defined by
\[
G(U):=\min\bigl\{m \in \mathbb{Z}_{\geq 0} \mid \exists \, U_1, \ldots, U_m \in E,\, U_m \cdots U_1 = U\bigr\},
\]
with the empty product understood as $\mathrm{Id}$ when $m = 0$. For $\epsilon > 0$, the \emph{approximate gate complexity} of $U$ with error $\epsilon$ is defined by
\[
G(U,\epsilon) := \min \bigl\{ m \in \mathbb{Z}_{\geq 0} \mid \exists\, U_1, \ldots, U_m \in E,\, \Vert U_m \cdots U_1 - U \Vert_{\mathrm{op}} < \epsilon \bigr\}.
\]
Note that every $U \in \mathrm{U}(D)$ is a finite product of elements of $E$~\cite[Section~3]{Knill1995}, and hence $G(U)$ and $G(U,\epsilon)$ are finite and well-defined.

	For the cliff metric $g^{\mathrm{SU}}_q$ on $\mathrm{SU}(D)$ with penalty $q$ chosen sufficiently large in terms of $N$ and $\epsilon$, Nielsen et al.~\cite{NielsenEtAl2006Geometry, Nielsen2006OptimalControl} established the polynomial equivalence
\[
d_q^{\mathrm{SU}}(\mathrm{Id}, U) \leq \pi G(U), \qquad G(U,\epsilon) = O\left(\frac{N^6d_q^{\mathrm{SU}}(\mathrm{Id}, U)^3}{\epsilon^2}\right).
\]
Here the first inequality holds on all of $\mathrm{U}(D)$ and for every $q > 1$ (as proved in Lemma~\ref{lem:diameter-upper-bound} below).

	We define the \emph{maximal exact gate complexity} and \emph{maximal approximate gate complexity} on $N$ qubits by
\[
G_{\max}(N) :=\sup_{U \in \mathrm{U}(D)} G(U), \qquad G_{\max}(N,\epsilon) := \sup_{U \in \mathrm{U}(D)} G(U, \epsilon).
\]
Knill~\cite[Theorems~3.1, 3.3, and~3.4]{Knill1995} proved that there exist constants $c_1, c_2 > 0$, independent of $N$, such that
\begin{equation}\label{equ:knill-bound}
c_14^N \leq G_{\max}(N) \leq c_2N4^N.
\end{equation}
Thus
\[
\lim_{N \to \infty} \frac{1}{N}\log_2 G_{\max}(N) = 2.
\]
His results also show that for every $\epsilon \in (0,\sqrt{2})$, the maximal approximate gate complexity satisfies~\cite[Corollary~4.2]{Knill1995}
\[
\lim_{N \to \infty} \frac{1}{N}\log_2 G_{\max}(N,\epsilon) = 2.
\]

	Combining the gate-complexity upper bound~\eqref{equ:knill-bound} with~\eqref{equ:diameter_upper_bound_1_of_cliff_metric}, we obtain the following diameter upper bounds for the finite- and infinite-cliff metrics.
	
\begin{lemma}\label{lem:diameter-upper-bound}
Let $N \geq 3$ be an integer, and set $D = 2^N$. For every $U \in \mathrm{U}(D)$ and every $q > 1$,
\begin{equation}\label{equ:distance-gate-upper-bound}
d_q(\mathrm{Id}, U) \leq d_\infty(\mathrm{Id}, U) \leq \pi G(U).
\end{equation}
Moreover, there exists a constant $C > 0$, independent of $N$ and $q$, such that
\begin{equation}\label{equ:diameter-upper-bound}
\mathrm{diam}(\mathrm{U}(D), d_\infty) \leq C N4^N, \qquad \mathrm{diam}(\mathrm{U}(D), d_q) \leq \min\bigl\{\pi\sqrt{q}, C N4^N\bigr\}.
\end{equation}
\end{lemma}

\begin{proof}
We first show that $d_\infty(\mathrm{Id}, U_0) \leq \pi$ for every $U_0 \in E$. By definition, there exist $S \subseteq \{1, \ldots, N\}$ with $|S| = 2$ and $W \in \mathrm{U}(\mathcal{V}_S)$ such that $U_0 = W \otimes \mathrm{Id}_{\mathcal{V}_{\overline{S}}}$ after reordering. Consider the Lie group embedding
\[
\iota_S: \mathrm{U}(\mathcal{V}_S) \hookrightarrow \mathrm{U}(D), \qquad \iota_S(A) := A\otimes \mathrm{Id}_{\mathcal{V}_{\overline{S}}}.
\]
Identifying $\mathcal{V}_S \cong \mathbb{C}^4$, its differential at the identity maps the $g_1$-orthonormal basis 
\[
\{-\sqrt{-1}\sigma_{j_1} \otimes \sigma_{j_2} \mid j_1,j_2 \in \{0,1,2,3\}\}
\]
of $\mathfrak{u}(4)$ into the $g_\infty$-orthonormal basis $\{\tilde\sigma_I\}_{I \in \mathcal{P}_N}$ of $\mathcal{H}_N$. By right-invariance, the tangent space of $\iota_S(\mathrm{U}(\mathcal{V}_S))$ lies in $\Delta_{\mathcal{H}_N}$ and $\iota_S^* g_\infty$ coincides with the normalized bi-invariant Riemannian metric $g_1$ on $\mathrm{U}(4)$. Hence $\iota_S$ pushes smooth curves in $\mathrm{U}(\mathcal{V}_S)$ forward to admissible curves of the same length in $(\mathrm{U}(D), \Delta_{\mathcal{H}_N}, g_\infty)$, and~\eqref{equ:diameter_of_bi_invariant_metric} gives
\[
d_\infty(\mathrm{Id}, U_0) \leq d_1(\mathrm{Id}_{\mathcal{V}_S}, W) \leq \mathrm{diam}(\mathrm{U}(4), d_1) = \pi.
\]

	Now let $U \in \mathrm{U}(D)$ and set $m := G(U)$. Write $U = U_m \cdots U_1$ with $U_1, \ldots, U_m \in E$. Then by the triangle inequality and the right-invariance of $d_\infty$,
\[
d_\infty(\mathrm{Id}, U) \leq \sum_{j = 1}^m d_\infty(U_j \cdots U_1, U_{j - 1} \cdots U_1) = \sum_{j = 1}^m d_\infty(\mathrm{Id}, U_j) \leq \pi G(U).
\]
Together with~\eqref{equ:metric_comparison}, this proves~\eqref{equ:distance-gate-upper-bound}.

	Finally, by the right-invariance of $d_\infty$, \eqref{equ:distance-gate-upper-bound}, and~\eqref{equ:knill-bound}, we obtain
\[
\mathrm{diam}(\mathrm{U}(D), d_\infty) = \sup_{U \in \mathrm{U}(D)} d_\infty(\mathrm{Id}, U) \leq \pi G_{\max}(N) \leq C N4^N
\]
with $C:= c_2\pi$ independent of $N$ and $q$. Combining this bound with~\eqref{equ:metric_comparison} and~\eqref{equ:diameter_upper_bound_1_of_cliff_metric} yields $\mathrm{diam}(\mathrm{U}(D), d_q) \leq \min\{\pi\sqrt{q}, C N4^N\}$, completing the proof.
\end{proof}

\subsection{Dyson expansion and unitary approximation}

	In this subsection, we introduce the truncated Dyson expansion~\cite{Dyson1949}, which is also known as the Volterra expansion in geometric control theory~\cite[Chapter~6]{AgrachevBarilariBoscain2020}, and establish three approximation estimates used in the covering argument of Section~3.	
	
	Fix an integer $N \geq 3$, set $D = 2^N$, and let $u = (u_I)_{I \in \mathcal{P}_N} \in L^\infty([0,1];\mathbb{R}^{b_N})$ be a control. Modifying each component $u_I$ on a set of measure zero, we may assume that $|u_I(t)| \leq \Vert u_I\Vert_{L^\infty([0,1])}$ for every $t \in [0,1]$ and every $I \in \mathcal{P}_N$. For $t \in [0,1]$, define
\[
X(t) := \sum_{I \in \mathcal{P}_N} u_I(t) \tilde \sigma_I \in \mathcal{H}_N, \qquad \Vert u \Vert_{1,t} := \int_0^t \sum_{I \in \mathcal{P}_N} |u_I(\tau)| \dd \tau,
\]
and consider the initial value problem
\begin{equation}\label{equ:right-invariant-system}
\begin{aligned}
\dot \gamma(t) &= X(t)\gamma(t), \qquad \text{for almost every }t \in [0,1],\\ \gamma(0) &= \mathrm{Id}.
\end{aligned}
\end{equation}
Since $X(t) \in \mathfrak{u}(D)$ for every $t \in [0,1]$, any solution of~\eqref{equ:right-invariant-system} in $M_D(\mathbb{C})$ remains in the compact group $\mathrm{U}(D)$. Hence, by~\cite[Appendix~C, Theorem~54 and Proposition~C.3.6]{Sontag1998}, \eqref{equ:right-invariant-system} admits a unique solution $\gamma: [0,1] \to \mathrm{U}(D)$.

	Moreover, since $\Vert X(t)\gamma(t) \Vert_{\mathrm{op}} = \Vert X(t) \Vert_{\mathrm{op}} \leq \sum_{I \in \mathcal{P}_N} \Vert u_I \Vert_{L^\infty([0,1])}$ for every $t \in [0,1]$ by~\eqref{equ:unitary_Pauli_matrices}, Lemmas~\ref{lem:basic_property_of_operator_norm} and~\ref{lem:compare-distance} yield
\[
d_1(\gamma(\tau),\gamma(\tau_1)) \leq \frac{\pi}{2}\Vert \gamma(\tau_1) - \gamma(\tau) \Vert_{\mathrm{op}} \leq \frac{\pi}{2}\left(\sum_{I \in \mathcal{P}_N} \Vert u_I \Vert_{L^\infty([0,1])}\right)|\tau_1 - \tau|, \qquad \tau, \tau_1 \in [0,1],
\]
so $\gamma$ is $d_1$-Lipschitz and hence admissible for the infinite-cliff metric.

	For every $t \in [0,1]$, $\gamma$ satisfies the integral equation $\gamma(t) = \mathrm{Id} + \int_0^t X(t_1)\gamma(t_1) \dd t_1$. Iterating this identity $k$ times for any integer $k \geq 0$, we obtain
\[
\begin{aligned}
\gamma(t) &= \mathrm{Id} + \int_0^t X(t_1)\gamma(t_1) \dd t_1 \\&= \mathrm{Id} + \int_0^t X(t_1) \left(\mathrm{Id} + \int_0^{t_1} X(t_2)\gamma(t_2) \dd t_2\right)\dd t_1 \\&= \cdots \\&=\mathrm{Id} + \sum_{n = 1}^k \int_0^{t} \dd t_1 \cdots \int_0^{t_{n - 1}} \dd t_n \,X(t_1) \cdots X(t_n) + R_{k + 1}(t),
\end{aligned}
\]
where the remainder is given by
\[
R_{k + 1}(t) := \int_0^t \dd t_1 \cdots \int_0^{t_k} \dd t_{k + 1} \, X(t_1) \cdots X(t_{k + 1})\gamma(t_{k + 1}).
\]
Using $X(t_i) = \sum_{I \in \mathcal{P}_N} u_{I}(t_i) \tilde\sigma_{I}$, we can rewrite this as the $k$-th order \emph{truncated Dyson expansion}
\begin{equation}\label{equ:dyson-expansion}
\gamma(t) = \mathrm{Id} + \sum_{n = 1}^k \sum_{\boldsymbol{\alpha} \in \mathcal{P}_N^n} a_{\boldsymbol{\alpha}}(t) T_{\boldsymbol{\alpha}} + R_{k + 1}(t),
\end{equation}
where for each integer $n \geq 1$ and each multi-index $\boldsymbol{\alpha} = (I_1, \ldots, I_n) \in \mathcal{P}_N^n$ (writing $|\boldsymbol{\alpha}| := n$),
\[
a_{\boldsymbol{\alpha}}(t) := \int_0^t \dd t_1 \int_0^{t_1} \dd t_2 \cdots \int_0^{t_{n - 1}} \dd t_n \, u_{I_1}(t_1) \cdots u_{I_n}(t_n), \qquad T_{\boldsymbol{\alpha}} := \tilde\sigma_{I_1} \cdots \tilde\sigma_{I_n}.
\]
Since each $\tilde\sigma_I$ is unitary, each $T_{\boldsymbol{\alpha}}$ belongs to $\mathrm{U}(D)$ and hence satisfies $\Vert T_{\boldsymbol{\alpha}} \Vert_{\mathrm{op}} = 1$ by Lemma~\ref{lem:basic_property_of_operator_norm}.

\begin{lemma}\label{lem:dyson-expansion}
	With the notation above, for any control $u \in L^\infty([0,1]; \mathbb{R}^{b_N})$, all integers $n \geq 1$ and $k \geq 0$, and all $t \in [0,1]$, the coefficients $a_{\boldsymbol{\alpha}}(t)$ and the remainder $R_{k + 1}(t)$ in the truncated Dyson expansion~\eqref{equ:dyson-expansion} satisfy
\begin{equation}\label{equ:Dyson-expansion-estimate}
\sum_{\boldsymbol{\alpha} \in \mathcal{P}_N^n} |a_{\boldsymbol{\alpha}}(t)| \leq \frac{\Vert u \Vert_{1,t}^n}{n!}, \qquad \Vert R_{k + 1}(t) \Vert_{\mathrm{op}} \leq \frac{\Vert u \Vert_{1,t}^{k + 1}}{(k + 1)!}.
\end{equation}
\end{lemma}

\begin{proof}

	Set $f(\tau) := \sum_{I \in \mathcal{P}_N} |u_I(\tau)|$, so that $\int_0^t f(\tau) \dd \tau = \Vert u \Vert_{1,t}$. For any integer $n \geq 1$ and $t \in [0,1]$, let $\Delta_{n,\rho}(t) := \{(t_1, \ldots, t_n) \in [0,t]^n \mid t > t_{\rho(1)} > \cdots > t_{\rho(n)} > 0\}$ for each permutation $\rho \in S_n$ of $\{1,\ldots,n\}$. The $n!$ sets $\{\Delta_{n,\rho}(t)\}_{\rho \in S_n}$ are pairwise disjoint, and their complement in $[0,t]^n$ consists of the boundary of $[0,t]^n$ together with the hyperplanes $\{t_i = t_j\}$ $(i \neq j)$, which has Lebesgue measure zero. Since $\prod_{i = 1}^n f(t_i)$ is invariant under coordinate permutations, we obtain
\begin{equation}\label{equ:simplex-integral}
\begin{aligned}
\int_0^t \dd t_1 \int_0^{t_1} \dd t_2 \cdots \int_0^{t_{n - 1}} \dd t_n \prod_{i = 1}^n f(t_i) &= \frac{1}{n!}\sum_{\rho \in S_n} \int_{\Delta_{n, \rho}(t)} \prod_{i = 1}^n f(t_i) \dd t_1 \cdots \dd t_n \\
&= \frac{1}{n!} \int_{[0,t]^n} \prod_{i = 1}^n f(t_i) \dd t_1 \cdots \dd t_n = \frac{\Vert u\Vert_{1,t}^n}{n!}.
\end{aligned}
\end{equation}

		Using~\eqref{equ:simplex-integral}, we obtain
\[
\begin{aligned}
\sum_{\boldsymbol{\alpha} \in \mathcal{P}_N^n}|a_{\boldsymbol{\alpha}}(t)| &= \sum_{\boldsymbol{\alpha} \in \mathcal{P}_N^n}\left|\int_0^t \dd t_1 \int_0^{t_1} \dd t_2 \cdots \int_0^{t_{n - 1}} \dd t_n\, u_{I_1}(t_1) \cdots u_{I_n}(t_n)\right| \\
&\leq \int_0^{t} \dd t_1 \int_0^{t_1} \dd t_2 \cdots \int_0^{t_{n - 1}} \dd t_{n} \sum_{\boldsymbol{\alpha} \in \mathcal{P}_N^n}\left|u_{I_1}(t_1) \cdots u_{I_n}(t_n)\right| \\
&= \int_0^{t} \dd t_1 \int_0^{t_1} \dd t_2 \cdots \int_0^{t_{n - 1}} \dd t_{n}\, \prod_{i = 1}^n f(t_i) \\
&= \frac{\Vert u\Vert_{1,t}^n}{n!}.
\end{aligned}
\]

	For the remainder estimate, combining~\eqref{equ:unitary_Pauli_matrices} with Lemma~\ref{lem:basic_property_of_operator_norm} gives
\[
\Vert X(\tau) \Vert_{\mathrm{op}} = \left\Vert \sum_{I \in \mathcal{P}_N} u_I(\tau)\tilde\sigma_I \right\Vert_{\mathrm{op}} \leq \sum_{I \in \mathcal{P}_N} |u_I(\tau)|\,\Vert \tilde\sigma_I \Vert_{\mathrm{op}} = f(\tau)
\]
for every $\tau \in [0,1]$. Since $\gamma(t_{k+1}) \in \mathrm{U}(D)$, applying Lemma~\ref{lem:basic_property_of_operator_norm} and~\eqref{equ:simplex-integral} gives
\[
\begin{aligned}
\Vert R_{k + 1}(t) \Vert_{\mathrm{op}} &\leq \int_0^{t} \dd t_1 \cdots \int_0^{t_k} \dd t_{k + 1} \,\Vert X(t_1) \Vert_{\mathrm{op}} \cdots \Vert X(t_{k + 1})\Vert_{\mathrm{op}} \\
&\leq \int_0^{t} \dd t_1 \cdots \int_0^{t_k} \dd t_{k + 1} \,\prod_{i = 1}^{k + 1} f(t_i) \\
&= \frac{\Vert u \Vert_{1,t}^{k + 1}}{(k + 1)!}.
\end{aligned}
\]
\end{proof}

	Without the remainder term, a truncated Dyson expansion need not be unitary. The following estimate allows us to replace a matrix approximation to a unitary matrix by a unitary approximation.

\begin{lemma}\label{lem:SVD}
	Let $D \geq 2$ be an integer, let $V \in M_D(\mathbb{C})$, and fix a singular value decomposition (SVD)
\[
V = A \Sigma B^\dagger,
\]
where $A, B\in \mathrm{U}(D)$ and $\Sigma$ is diagonal with nonnegative entries. Set $Q = AB^\dagger$. Then $Q \in \mathrm{U}(D)$, and for every $F \in \mathrm{U}(D)$,
\[
\Vert F - Q\Vert_{\mathrm{op}} \leq 2 \Vert F - V \Vert_{\mathrm{op}}.
\]
\end{lemma}

\begin{proof}
	Set $H = B\Sigma B^\dagger$. Then $H$ is positive semidefinite, $Q$ is unitary, and $V = QH$. Since the operator norm is unitarily invariant by Lemma~\ref{lem:basic_property_of_operator_norm}, \cite[Theorem~1]{FanHoffman1955} gives
\[
\Vert V - Q\Vert_{\mathrm{op}} \leq \Vert V - F \Vert_{\mathrm{op}}, \qquad \text{for every }F \in \mathrm{U}(D).
\]
The triangle inequality therefore yields
\[
\Vert F - Q\Vert_{\mathrm{op}} \leq \Vert F - V \Vert_{\mathrm{op}} + \Vert V - Q\Vert_{\mathrm{op}} \leq 2\Vert F - V\Vert_{\mathrm{op}}.
\]
\end{proof}

	The next estimate controls how approximation errors accumulate when unitary matrices are multiplied.
	
\begin{lemma}
\label{lem:telescope}
Let $D \geq 2$ and $S \geq 1$ be integers, and let $F_1, \ldots, F_S, Q_1, \ldots, Q_S \in \mathrm{U}(D)$. Then
\[
\Vert F_S \cdots F_1 - Q_S \cdots Q_1 \Vert_{\mathrm{op}} \leq \sum_{j = 1}^S \Vert F_j - Q_j \Vert_{\mathrm{op}}.
\]
\end{lemma}
\begin{proof}
Note that
\[
F_S \cdots F_1 - Q_S \cdots Q_1 = \sum_{j = 1}^S F_S \cdots F_{j+1}(F_j - Q_j)Q_{j - 1} \cdots Q_1.
\]
By Lemma~\ref{lem:basic_property_of_operator_norm},
\[
\Vert F_S \cdots F_1 - Q_S \cdots Q_1 \Vert_{\mathrm{op}} \leq \sum_{j = 1}^S \Vert F_S \cdots F_{j+1}(F_j - Q_j)Q_{j - 1} \cdots Q_1 \Vert_{\mathrm{op}} = \sum_{j = 1}^S \Vert F_{j} - Q_{j}\Vert_{\mathrm{op}}.
\]
\end{proof}

\section{Covering estimates for finite-cliff balls}

	In this section, we establish upper bounds on the $d_1$-covering numbers of finite-cliff balls $B_q(\mathrm{Id}, L)$. The basic intuition is the large gap between $b_N = \dim \mathcal{H}_N = O(N^2)$ and the dimension of the unitary group $\dim \mathrm{U}(2^N) = 4^N$, which allows $B_q(\mathrm{Id}, L)$ to be approximated by a relatively small finite net in $\mathrm{U}(2^N)$.
	
	A core step is Lemma~\ref{lem:distance-upper-bound}, which constructs a finite net $\mathcal{C}_{S,K,k} \subseteq \mathrm{U}(2^N)$ approximating $B_q(\mathrm{Id}, L)$ in $d_1$. Intuitively, we first use horizontal projection to reduce from $B_q(\mathrm{Id}, L)$ to the infinite-cliff ball $B_\infty(\mathrm{Id}, L)$ (with a $d_1$-error less than $L/\sqrt{q}$), and then discretize the projected admissible curve via the truncated Dyson expansion. Here $S$ is the number of equal-length sub-arcs into which the projected curve is subdivided, $k$ is the truncation order of the Dyson expansion on each sub-arc, and $2K$ is the number of equal subintervals used to discretize the interval containing each expansion coefficient.

\begin{lemma}\label{lem:distance-upper-bound}
	Let $N \geq 3$ be an integer and set $D = 2^N$. Equip $\mathrm{U}(D)$ with the cliff metric $g_q$, $q > 1$. Fix $L > 0$. For any integers $S,K,k \geq 1$, set $h = \frac{\sqrt{b_N}L}{S}$. Then there exists a finite subset $\mathcal{C}_{S,K,k} \subseteq \mathrm{U}(D)$ with $|\mathcal{C}_{S,K,k}| \leq (2K + 1)^{S\sum_{n = 1}^k b_N^n}$ such that, for every $x \in B_q(\mathrm{Id},L)$, there exists $y \in \mathcal{C}_{S,K,k}$ satisfying
\begin{equation}\label{equ:distance-upper-bound}
d_1(x,y) < \frac{L}{\sqrt{q}} + \frac{\pi}{2}\left(\frac{S}{K}\sum_{n = 1}^k \frac{(b_Nh)^n}{n!} + \frac{2Sh^{k + 1}}{(k + 1)!}\right).
\end{equation}
\end{lemma}

\begin{proof}
	Fix $x \in B_q(\mathrm{Id}, L)$ and choose a smooth curve $\gamma: [0,1] \to \mathrm{U}(D)$ joining $\mathrm{Id}$ to $x$ with
\[
l_{g_q}(\gamma) < L.
\]
Write $X(t) = \dot \gamma(t)\gamma(t)^\dagger$, and let $\gamma_P$ be the solution of
\[
\dot\gamma_P(t) = P_N(X(t)) \gamma_P(t), \qquad \gamma_P(0) = \mathrm{Id}.
\]
By the definition of $g_q$,
\[
\Vert X(t) \Vert_{g_q}^2 = \Vert P_N(X(t)) \Vert_{g_1}^2 + q \Vert X(t) - P_N(X(t)) \Vert_{g_1}^2.
\]
Since $P_N(X(t)) \in \mathcal{H}_N$, the curve $\gamma_P$ is admissible for the infinite-cliff metric, and
\[
l_{\infty}(\gamma_P) = \int_0^1\Vert P_N(X(t)) \Vert_{g_1} \dd t \leq \int_0^1 \Vert X(t)\Vert_{g_q} \dd t = l_{g_q}(\gamma) < L.
\]
Thus $\tilde x := \gamma_P(1)$ belongs to $B_\infty(\mathrm{Id}, L)$.

	Set $w(t) = \gamma_P(t)^\dagger \gamma(t)$. Then $w(0) = \mathrm{Id}$ and $\dot w = \gamma_P^\dagger (X - P_N(X))\gamma_P w$. By the bi-invariance of $g_1$, we obtain
\[
\begin{aligned}
d_1(x,\tilde x) &= d_1(\mathrm{Id},w(1)) \leq l_{g_1}(w) \\
&= \int_0^1 \Vert X(t) - P_N(X(t)) \Vert_{g_1} \dd t \\
&\leq \frac{1}{\sqrt{q}} \int_0^1 \Vert X(t) \Vert_{g_q} \dd t < \frac{L}{\sqrt{q}}.
\end{aligned}
\]

	For $1 \leq n \leq k$, define
\[
\begin{aligned}
\Lambda^{(n)}_{S,K} &:= \left\{\frac{h^n}{n!}\frac{j}{K} \middle| -K\leq j \leq K, j \in \mathbb{Z}\right\};\\
\Phi_{S,K}^{(n)} &:= \{\phi^{(n)}: \mathcal{P}_N^n \to \Lambda_{S,K}^{(n)}\};\\
\Phi_{S,K,k} &:= \prod_{n = 1}^k \Phi_{S,K}^{(n)}.
\end{aligned}
\]
We write an element $\phi \in \Phi_{S,K,k}$ as $\phi = (\phi^{(1)}, \ldots, \phi^{(k)})$, where $\phi^{(n)} \in \Phi_{S,K}^{(n)}$. For $1 \leq n \leq k$ and $\boldsymbol{\alpha} \in \mathcal{P}_N^n$, define $\phi(\boldsymbol{\alpha}) := \phi^{(n)}(\boldsymbol{\alpha})$.

	For $\phi \in \Phi_{S,K,k}$, we can construct a unitary matrix $Q_\phi$ as follows. Define $F_{\phi}$ by
\[
F_\phi := \mathrm{Id} + \sum_{n = 1}^{k} \sum_{\boldsymbol{\alpha} \in \mathcal{P}_N^n} \phi(\boldsymbol{\alpha}) T_{\boldsymbol{\alpha}}.
\]
	Fix an SVD $F_{\phi} = A_{\phi}\Sigma_{\phi}B_{\phi}^\dagger$ and define $Q_{\phi} = A_{\phi}B^\dagger_{\phi}$. Furthermore, define
\[
\mathcal{C}_{S,K,k} :=\{Q_{\phi_S} \cdots Q_{\phi_1} \mid \phi_j \in \Phi_{S,K,k} \text{ for } j = 1, \ldots, S\}.
\]
Here
\[
|\mathcal{C}_{S,K,k}| \leq |\Phi_{S,K,k}|^S = (2K+1)^{S\sum_{n = 1}^k b_N^n}.
\]
	We next approximate $\tilde x$ by an element of $\mathcal{C}_{S,K,k}$ in the operator norm. Let $\eta: [0,1] \to \mathrm{U}(D)$ be a constant-speed reparametrization of $\gamma_P$, and set
\[
\ell := l_\infty(\eta) = l_{\infty}(\gamma_P) < L.
\]
Let $u = (u_I)_{I \in \mathcal{P}_N}$ be its control relative to $\{\tilde\sigma_I\}_{I \in \mathcal{P}_N}$. Then, for almost every $t \in [0,1]$,
\[
\dot \eta(t) = \left(\sum_{I \in \mathcal{P}_N} u_{I}(t)\tilde \sigma_{I}\right)\eta(t), \qquad \left(\sum_{I \in \mathcal{P}_N}u_I(t)^2 \right)^{1/2} = \ell.
\]

	For $j = 1, \ldots, S$, define
\[
z_j(t) := \eta\left(\frac{t + (j - 1)}{S}\right)\eta\left(\frac{j - 1}{S}\right)^\dagger, \qquad t \in [0,1].
\]
Then $z_j(0) = \mathrm{Id}$, and its control is
\[
u_{I,j}(t) = \frac{1}{S}u_{I} (\frac{t + (j - 1)}{S}).
\]
	Since $\sum_{I \in \mathcal{P}_N} |u_I(\tau)| \leq \sqrt{b_N}\bigl(\sum_{I \in \mathcal{P}_N} u_I(\tau)^2\bigr)^{1/2} = \sqrt{b_N}\ell$ for almost every $\tau \in [0,1]$, we obtain
\begin{equation}\label{equ:estimate-u_j}
\begin{aligned}
\Vert u_j \Vert_{1,1} &= \int_{0}^1 \sum_{I \in \mathcal{P}_N} |u_{I,j}(t)|\dd t \\& = \int_{0}^1 \frac{1}{S} \sum_{I \in \mathcal{P}_N} \left|u_{I}(\frac{t + (j - 1)}{S})\right|\dd t \\&\leq \frac{1}{S} \int_0^1 \sqrt{b_N}\ell \dd t \\&< \frac{\sqrt{b_N}L}{S} = h.
\end{aligned}
\end{equation}
Applying the Dyson expansion to $z_j$ at $t = 1$, we obtain
\[
z_j(1) = \mathrm{Id} + \sum_{n = 1}^k \sum_{\boldsymbol{\alpha} \in \mathcal{P}_N^n} a_{\boldsymbol{\alpha},j}(1) T_{\boldsymbol{\alpha}} + R_{k+1,j}(1).
\]

	Lemma~\ref{lem:dyson-expansion} gives
\begin{equation}\label{equ:estimate-dyson-set}
\sum_{\boldsymbol{\alpha} \in \mathcal{P}_N^n}|a_{\boldsymbol{\alpha},j}(1)| \leq \frac{\Vert u_{j} \Vert_{1,1}^{n}}{n!} < \frac{h^{n}}{n!}, \qquad \Vert R_{k+1,j}(1) \Vert_{\mathrm{op}} \leq \frac{\Vert u_j \Vert_{1,1}^{k + 1}}{(k + 1)!} < \frac{h^{k + 1}}{(k + 1)!}.
\end{equation}
Thus we can find $\phi_j \in \Phi_{S,K,k}$ such that for all $1 \leq n \leq k$, $\boldsymbol{\alpha} \in \mathcal{P}_N^n$,
\begin{equation}\label{equ:choice-of-point}
|\phi_j(\boldsymbol{\alpha}) - a_{\boldsymbol{\alpha},j}(1)| \leq \frac{1}{2K}\frac{h^{n}}{n!}.
\end{equation}

	Combining \eqref{equ:estimate-dyson-set} and \eqref{equ:choice-of-point} for $Q_{\phi_j}$ and applying Lemma~\ref{lem:SVD},  we obtain
\[
\begin{aligned}
\Vert Q_{\phi_j} - z_j(1) \Vert_{\mathrm{op}} &\leq 2\Vert F_{\phi_j} - z_j(1) \Vert_{\mathrm{op}} \\
&\leq 2\left(\sum_{n = 1}^k \sum_{\boldsymbol{\alpha} \in \mathcal{P}_N^n} |\phi_j(\boldsymbol{\alpha})-a_{\boldsymbol{\alpha},j}(1)| \Vert T_{\boldsymbol{\alpha}} \Vert_{\mathrm{op}} + \Vert R_{k + 1,j}(1) \Vert_{\mathrm{op}}\right) \\
&< 2\left(\sum_{n = 1}^k b_N^n \cdot \frac{1}{2K}\frac{h^n}{n!} + \frac{h^{k + 1}}{(k + 1)!}\right) \\
&= \frac{1}{K}\sum_{n = 1}^k\frac{(b_Nh)^n}{n!} + \frac{2h^{k + 1}}{(k + 1)!}.
\end{aligned}
\]
Note that $\tilde x = z_S(1) \cdots z_1(1)$. Let $y := Q_{\phi_S} \cdots Q_{\phi_1} \in \mathcal{C}_{S,K,k}$. Thus by Lemma~\ref{lem:telescope}, we have
\[
\Vert \tilde x - y \Vert_{\mathrm{op}} \leq \sum_{j = 1}^S \Vert z_j(1) - Q_{\phi_j} \Vert_{\mathrm{op}} < \frac{S}{K}\sum_{n = 1}^k \frac{(b_Nh)^n}{n!} + \frac{2Sh^{k + 1}}{(k + 1)!}.
\]
Finally, Lemma~\ref{lem:compare-distance} gives
\[
\begin{aligned}
d_1(x,y) &\leq d_1(x, \tilde x) + d_1(\tilde x,y) \\&\leq d_1(x, \tilde x) + \frac{\pi}{2}\Vert \tilde x - y \Vert_{\mathrm{op}} \\&< \frac{L}{\sqrt{q}} + \frac{\pi}{2}\left(\frac{S}{K}\sum_{n = 1}^k \frac{(b_Nh)^n}{n!} + \frac{2Sh^{k + 1}}{(k + 1)!}\right).
\end{aligned}
\]
\end{proof}

We now choose the parameters $S$ and $K$ in Lemma~\ref{lem:distance-upper-bound} to obtain the following bound on the $d_1$-covering numbers of finite-cliff balls.

\begin{proposition}\label{prop:upper-bound-of-covering-number}
	For every integer $k \geq 1$ and every $\epsilon > 0$, there exists a constant $C_{k,\epsilon} > 0$, depending only on $k$ and $\epsilon$, such that, for all integers $N \geq 3$, all $q > 1$, and all $L > 0$,
\begin{equation}\label{equ:upper-bound-of-covering-number}
\ln (N_{1}(B_q(\mathrm{Id},L), \frac{\pi}{2}\epsilon + \frac{L}{\sqrt{q}})) \leq C_{k,\epsilon} (b_N^k + b_N^{k + 3/2}L + b_N^{k + 1/2 + \frac{1}{2k}}L^{1 + \frac{1}{k}})\ln(3 + \frac{4eb_N^{3/2}L}{\epsilon}).
\end{equation}
\end{proposition}

\begin{proof}

	We apply Lemma~\ref{lem:distance-upper-bound}. For integers $S,K \geq 1$, we can find $\mathcal{C}_{S,K,k}$ such that for every $x \in B_{q}(\mathrm{Id},L)$, there exists $y \in \mathcal{C}_{S,K,k}$ satisfying
\[
d_1(x,y)  < \frac{L}{\sqrt{q}} + \frac{\pi}{2}\left( \frac{S}{K}\sum_{n = 1}^k \frac{(b_Nh)^n}{n!} + \frac{2Sh^{k + 1}}{(k + 1)!}\right),
\]
where $h = \frac{\sqrt{b_N}L}{S}$.

	Let $h_{0,k} = \min\left\{\frac{1}{b_N}, \left(\frac{\epsilon(k + 1)!}{4\sqrt{b_N}L}\right)^{1/k}\right\}$ and $S = \left\lceil \frac{\sqrt{b_N}L}{h_{0,k}} \right\rceil$. Then
\[
h = \frac{\sqrt{b_N}L}{S} \leq h_{0,k} \leq \frac{1}{b_N},
\]
and hence $0 < b_N h \leq 1$. Then
\[
2S \frac{h^{k + 1}}{(k + 1)!} \leq 2Sh \frac{h_{0,k}^k}{(k + 1)!} \leq 2\sqrt{b_N}L \frac{\epsilon}{4\sqrt{b_N}L} = \frac{\epsilon}{2}.
\]
Choose $K = \left\lceil \frac{2e b_N^{3/2}L}{\epsilon} \right\rceil$. Since $0 < b_Nh \leq 1$, we have
\[
\frac{S}{K} \sum_{n = 1}^k \frac{(b_Nh)^n}{n!} \leq \frac{S}{K}(e^{b_Nh} - 1) < \frac{Sehb_N}{K} = \frac{eb_N^{3/2}L}{K} \leq \frac{\epsilon}{2}.
\]
Thus we can see that
\[
\frac{S}{K}\sum_{n = 1}^k \frac{(b_Nh)^n}{n!} + \frac{2Sh^{k + 1}}{(k + 1)!} < \epsilon.
\]
Hence
\[
\begin{aligned}
N_{1}\left(B_{q}(\mathrm{Id},L),\frac{\pi}{2}\epsilon + \frac{L}{\sqrt{q}}\right) \leq |\mathcal{C}_{S,K,k}| \leq (2K + 1)^{S\sum_{n = 1}^k b_N^n},
\end{aligned}
\]
where $S \leq1 + \frac{\sqrt{b_N}L}{h_{0,k}} \leq 1 + b_N^{3/2}L + \bigl(\frac{4}{\epsilon(k + 1)!}\bigr)^{1/k}b_N^{1/2 + \frac{1}{2k}}L^{1 + \frac{1}{k}}$,
$2K + 1 \leq 3 + \frac{4eb_N^{3/2}L}{\epsilon}$, and
$\sum_{n = 1}^k b_N^n \leq kb_N^k$. Setting $C_{k,\epsilon} := k + k\left(\frac{4}{\epsilon(k + 1)!}\right)^{1/k}$, we obtain
\begin{align*}
&\ln (N_{1}(B_q(\mathrm{Id},L),\frac{\pi}{2}\epsilon + \frac{L}{\sqrt{q}})) \\&\leq S\left(\sum_{n = 1}^k b_N^n\right)\ln(2K + 1) \\&\leq C_{k,\epsilon} (b_N^k + b_N^{k + 3/2}L + b_N^{k + 1/2 + \frac{1}{2k}}L^{1 + \frac{1}{k}})\ln(3 + \frac{4eb_N^{3/2}L}{\epsilon}).\qedhere
\end{align*}
\end{proof}

\section{Proof of the main results}

	We now combine the covering estimates from Section~3 with the Haar-measure bound for $d_1$-balls to establish the diameter lower bounds and the corresponding measure estimates for the cliff metrics. We then combine these estimates with the diameter upper bounds from Lemma~\ref{lem:diameter-upper-bound} to complete the proof of Theorem~\ref{thm:main-theorem-intro}.

\begin{proposition}\label{prop:lower-bound-finite-cliff}
	Fix $\lambda > 0$. For each integer $N \geq 3$, let $q_{N,\lambda} = 2^{\lambda N}$. Then
\begin{equation}\label{equ:lower-bound-finite-cliff}
\liminf_{N \to \infty}\frac{1}{N} \log_2 \mathrm{diam}(\mathrm{U}(2^N), d_{q_{N,\lambda}}) \geq \min\{\frac{\lambda}{2},2\}.
\end{equation}
Furthermore, for every $\eta < \min\{\frac{\lambda}{2}, 2\}$,
\[
\mu_{2^N}(B_{q_{N,\lambda}}(\mathrm{Id}, 2^{\eta N})) \longrightarrow 0 ,\qquad \text{as } N \to \infty.
\]
\end{proposition}

\begin{proof}
	Note that for $\eta_1 \leq \eta_2$, we have $B_{q_{N,\lambda}}(\mathrm{Id}, 2^{\eta_1 N})\subseteq B_{q_{N,\lambda}}(\mathrm{Id}, 2^{\eta_2 N})$. Since $\min\{\frac{\lambda}{2},2\} > 0$, it suffices to prove the Haar-measure assertion for $0 < \eta < \min\{\frac{\lambda}{2},2\}$.

	Fix $0 < \eta < \min \{\frac{\lambda}{2}, 2\}$ and choose an integer $k \geq 1$ such that
\[
\eta (1 + \frac{1}{k}) < 2.
\]
Set $L_N = 2^{\eta N}$. Now fix $0 < r \leq \epsilon_0$, where $\epsilon_0$ is as in Lemma~\ref{lem:unitary-bi-inv-riemannian-small-ball}. There exists a constant $c_0 > 0$, independent of $N$, such that
\[
\mu_{2^N}(B_1(x,r)) \leq \exp(-c_02^{2N}),\qquad \text{for every }x \in \mathrm{U}(2^N).
\]
By the definition of the covering number $N_1$, we have
\[
\mu_{2^N}(B_{q_{N,\lambda}}(\mathrm{Id},L_N)) \leq N_1(B_{q_{N,\lambda}}(\mathrm{Id}, L_N), r) \exp(-c_0 2^{2N}).
\]

	It remains to estimate $N_1(B_{q_{N, \lambda}}(\mathrm{Id}, L_N), r)$. Set $\epsilon = \frac{r}{\pi}$. Since $\eta < \frac{\lambda}{2}$, we have
\[
\frac{L_N}{\sqrt{q_{N,\lambda}}} = 2^{(\eta - \lambda/2)N} \longrightarrow 0 ,\qquad \text{as } N \to \infty.
\]
Hence, for all sufficiently large $N$, we have $\frac{L_N}{\sqrt{q_{N,\lambda}}} \leq \frac{r}{2}$, and thus
\[
\frac{\pi}{2}\epsilon + \frac{L_N}{\sqrt{q_{N,\lambda}}} \leq \frac{r}{2} + \frac{r}{2} = r.
\]

	Therefore, by Proposition~\ref{prop:upper-bound-of-covering-number}, for all sufficiently large $N$, we have
\[
\begin{aligned}
\ln (N_1(B_{q_{N,\lambda}}(\mathrm{Id},L_N), r)) &\leq \ln (N_1(B_{q_{N,\lambda}}(\mathrm{Id}, L_N), \frac{\pi}{2}\epsilon + \frac{L_N}{\sqrt{q_{N,\lambda}}}))\\
&\leq C_{k,\epsilon}(b_N^k + b_N^{k + 3/2}2^{\eta N} + b_N^{k + 1/2 + \frac{1}{2k}}2^{\eta(1 + 1/k)N})\ln(3 + \frac{4e b_N^{3/2}2^{\eta N}}{\epsilon}).
\end{aligned}
\]

Since $b_N=O(N^2)$ and $k$, $\eta$ and $\epsilon$ are fixed, $\ln(3 + \frac{4e b_N^{3/2}2^{\eta N}}{\epsilon}) = O(N)$. The inequality $\eta(1 + \frac{1}{k}) < 2$ therefore implies
\[
\ln (N_1(B_{q_{N,\lambda}}(\mathrm{Id},L_N), r)) = o(2^{2N}).
\]
This yields
\[
\mu_{2^N}(B_{q_{N,\lambda}}(\mathrm{Id},L_N)) \leq \exp(-c_02^{2N} + o(2^{2N})) \longrightarrow 0 ,\qquad \text{as } N \to \infty.
\]
This proves the Haar-measure assertion.

	In particular, for all sufficiently large $N$, $\mu_{2^N}(B_{q_{N,\lambda}}(\mathrm{Id}, L_N)) < 1$, and hence \linebreak$\mathrm{diam}(\mathrm{U}(2^N), d_{q_{N,\lambda}}) \geq L_N$. Thus
\[
\liminf_{N \to \infty}\frac{1}{N}\log_2 \mathrm{diam}(\mathrm{U}(2^N), d_{q_{N,\lambda}}) \geq \eta.
\]
Since this holds for every $0 < \eta < \min\{\frac{\lambda}{2},2\}$, we conclude that
\[
\liminf_{N \to \infty}\frac{1}{N}\log_2\mathrm{diam}(\mathrm{U}(2^N), d_{q_{N,\lambda}}) \geq \min\{\frac{\lambda}{2},2\}.
\]
\end{proof}

\begin{proof}[Proof of Theorem~\ref{thm:main-theorem-intro}]

	We first prove~(1). The lower bound and the Haar-measure assertion in (1) follow from Proposition~\ref{prop:lower-bound-finite-cliff}. For the upper bound, applying~\eqref{equ:diameter-upper-bound} in Lemma~\ref{lem:diameter-upper-bound} with $q = q_{N,\lambda} = 2^{\lambda N}$ gives
\[
\mathrm{diam}(\mathrm{U}(2^N), d_{q_{N, \lambda}}) \leq \min \bigl\{\pi2^{\frac{\lambda}{2}N}, C N4^N \bigr\},
\]
where $C > 0$ is independent of $N$ and $\lambda$. Taking $\frac{1}{N}\log_2$ and letting $N \to \infty$, we obtain
\[
\limsup_{N \to \infty}\frac{1}{N}\log_2\mathrm{diam}(\mathrm{U}(2^N), d_{q_{N,\lambda}}) \leq \min\{\frac{\lambda}{2}, 2\}.
\]
Together with the lower bound~\eqref{equ:lower-bound-finite-cliff}, this shows
\[
\lim_{N \to \infty}\frac{1}{N}\log_2\mathrm{diam}(\mathrm{U}(2^N), d_{q_{N,\lambda}}) = \min\{\frac{\lambda}{2}, 2\}.
\]
	
	We now prove~(2). By~\eqref{equ:metric_comparison} and~\eqref{equ:diameter-upper-bound}, taking $q = q_{N,4} = 16^N$ gives
\[
\mathrm{diam}(\mathrm{U}(2^N), d_{q_{N,4}}) \leq \mathrm{diam}(\mathrm{U}(2^N), d_\infty) \leq C N4^N.
\]
Applying~(1) with $\lambda = 4$ therefore yields
\[
\lim_{N\to\infty}\frac{1}{N}\log_2\mathrm{diam}(\mathrm{U}(2^N),d_\infty)=2.
\]
Moreover, by~\eqref{equ:metric_comparison}, $B_\infty(\mathrm{Id}, 2^{\eta N}) \subseteq B_{q_{N,4}}(\mathrm{Id}, 2^{\eta N})$, and hence for every $\eta < 2$, part~(1) with $\lambda = 4$ implies
\[
\mu_{2^N}(B_\infty(\mathrm{Id},2^{\eta N})) \leq \mu_{2^N}(B_{q_{N,4}}(\mathrm{Id},2^{\eta N})) \longrightarrow 0, \qquad \text{as }N \to \infty.
\]
This completes the proof.
\end{proof}

\begin{ack}
Xiaobo Liu was partially supported by the National Natural Science Foundation of China under Grant Numbers 12341105 and 12526302.
\end{ack}

\vspace{30pt} \noindent
Xiaobo Liu \\
School of Mathematical Sciences \& \\
Beijing International Center for Mathematical Research, \\
Peking University, Beijing, China. \\
Email: {\it xbliu@math.pku.edu.cn}
\ \\ \ \\
Lei Zheng \\
School of Mathematical Sciences, \\
Peking University, Beijing, China. \\
Email: {\it zhenglei0813@stu.pku.edu.cn}
\end{document}